\documentclass[10pt]{amsart}

\usepackage{amssymb,amsmath,amsthm,amsfonts}
\usepackage[mathscr]{eucal}
\usepackage[all]{xy}

\theoremstyle{plain}
\newtheorem{Theorem}{Theorem}
\newtheorem{Proposition}[Theorem]{Proposition}

\newtheorem{Lemma}[Theorem]{Lemma}
\newtheorem{Corollary}[Theorem]{Corollary}

\theoremstyle{definition}
\newtheorem{Definition}[Theorem]{Definition}
\newtheorem{Example}[Theorem]{Example}
\newtheorem{Remark}[Theorem]{Remark}
\newtheorem{Problem}[Theorem]{Problem}
\newtheorem{Question}[Theorem]{Question}

\begin{document}
\title{Anti-tori in square complex groups}
\author{Diego Rattaggi}
\address{Universit\'e de Gen\`eve,
Section de math\'ematiques,
2--4 rue du Li\`evre, CP 64, 
CH--1211 Gen\`eve 4, Switzerland}
\email{rattaggi@math.unige.ch}
\renewcommand{\subjclassname}{\textup{2000} Mathematics Subject Classification}
\date{April 4, 2005}
\begin{abstract}
An anti-torus is a subgroup $\langle a,b \rangle$ in the fundamental group
of a compact non-positively curved space $X$, 
acting in a specific way on the universal covering space $\tilde{X}$
such that $a$ and $b$ do not have any commuting non-trivial powers. 
We construct and investigate anti-tori in a class of commutative transitive fundamental groups of 
finite square complexes, in particular for the groups
$\Gamma_{p,l}$ originally studied by Mozes \cite{Mozes1}. 
It turns out that anti-tori in $\Gamma_{p,l}$ directly correspond to
non-commuting pairs of Hamilton quaternions.
Moreover, free anti-tori in $\Gamma_{p,l}$ are related to free groups 
generated by two integer quaternions, and also to free subgroups of $\mathrm{SO}_3(\mathbb{Q})$.
As an application, we prove that the multiplicative group generated by
the two quaternions $1+2i$ and $1+4k$ is not free.
\end{abstract}
\maketitle

\setcounter{section}{-1}
\section{Introduction} \label{SectionIntro}
Bridson and Wise have given the following definition
of an anti-torus \cite[Definition 9.1]{BridsonWise}:
Let $X$ be a compact non-positively curved space with universal cover
$p: \tilde{X} \to X$.
It is well-known that the fundamental group $\pi_1 (X, x)$ acts on $\tilde{X}$,
and that each element $\gamma \in \pi_1 (X, x)$ leaves invariant in
this action
at least one isometrically embedded copy of the real line,
a so-called axis for $\gamma$.
Let $a, b \in \pi_1 (X, x)$ and
suppose that there is an isometrically embedded plane in $\tilde{X}$
which contains an axis for each of $a, b$ and that these axes
intersect in $p^{-1} x$.
If $a$ and $b$ do not have powers that commute, then 
$\langle a, b \rangle$ is called an \emph{anti-torus} in $\pi_1 (X, x)$.
If $\langle a, b \rangle$ is free then it is called 
a \emph{free anti-torus}.

We will restrict to a class where 
$\tilde{X} = \mathcal{T}_{2m} \times \mathcal{T}_{2n}$,
the product of two regular trees of degree $2m$ and $2n$,
respectively,
and $X$ is a certain finite square complex having a single 
vertex $x$. The fundamental group 
$\pi_1 (X, x) < \mathrm{Aut}(\mathcal{T}_{2m}) \times \mathrm{Aut}(\mathcal{T}_{2n})$
is then called a $(2m,2n)$--group 
(see Section~\ref{SectionPrel} for the precise definition).
 
Wise \cite{Wise} has constructed an anti-torus in a $(4,6)$--group to 
produce the first examples of non-residually finite groups in the following three important
classes: finitely presented small cancellation groups, automatic groups,
and groups acting properly discontinuously and cocompactly on CAT(0)-spaces.
Another application of anti-tori is the generation of
aperiodic tilings of the Euclidean plane by unit squares
(see \cite{Wise}, \cite{Rattaggi}).

In general, it seems to be very difficult to decide whether a subgroup
$\langle a, b \rangle$ is an anti-torus, or to decide
whether a group $\pi_1 (X, x)$ has an anti-torus or not.
In Section~\ref{SectionCT}, we further restrict to 
\emph{commutative transitive} $(2m,2n)$--groups, i.e.\ to groups 
$G$ where commutativity is a transitive relation on 
$G \setminus \{ 1 \}$.
In this context, we prove a dichotomy that 
$\langle a, b \rangle$ either is an anti-torus,
or isomorphic to the abelian group 
$\mathbb{Z} \times \mathbb{Z}$. Moreover, it turns out
that any commutative transitive 
$(2m,2n)$--group has an anti-torus, if $(m,n) \ne (1,1)$.
In Section~\ref{SectionGamma}, we define for any
pair $(p,l)$ of distinct odd prime numbers
a commutative transitive $(p+1,l+1)$--group $\Gamma_{p,l}$
and apply the results of Section~\ref{SectionCT}.
Anti-tori in $\Gamma_{p,l}$ are directly related to 
non-commuting Hamilton quaternions $x, y \in \mathbb{H}(\mathbb{Z})$
of norm a power of $p$ and $l$, respectively.
Although these considerations provide a very easy
method to construct anti-tori in $\Gamma_{p,l}$,
it is not clear at all if there are \emph{free} anti-tori
in $(2m,2n)$--groups. We give in Section~\ref{SectionFree}
a criterion for the construction of free anti-tori 
in terms of free groups generated by two quaternions,
but do not know if such quaternions exist.
Nevertheless, this criterion can be applied
to prove that certain pairs of quaternions, for example 
$1+2i$ and $1+4k$, do not generate a free group,
and we establish an explicit (long) relation
in this example.
Finally, we relate in Section~\ref{SectionSO}
free subgroups of $\Gamma_{p,l}$ to free subgroups
of $\mathrm{SO}_3(\mathbb{Q})$, using 
an explicit embedding $\Gamma_{p,l} \to \mathrm{SO}_3(\mathbb{Q})$.

Most results of this work are taken from the
authors Ph.D.\ thesis \cite{Rattaggi}.

\section{Preliminaries} \label{SectionPrel}
Let $m, n \in \mathbb{N}$ and $E_h := \{ a_1, \ldots , a_m \} ^{\pm 1}$,
$E_v := \{ b_1, \ldots , b_n \} ^{\pm 1}$.
A \emph{$(2m,2n)$--group} is the fundamental group $\Gamma = \pi_1 (X,x)$
of a finite $2$-dimensional cell complex $X$ satisfying the following conditions:
\begin{itemize}
\item The one-skeleton $X^{(1)}$ consists of a single vertex $x$ and
$m+n$ oriented loops $a_1, \ldots, a_m$, $b_1, \ldots, b_n$,
whose inverses are denoted by $a_1^{-1}, \ldots, a_m^{-1}$,
$b_1^{-1}, \ldots, b_n^{-1}$. In other words, $X^{(1)}$ is the graph 
with vertex set $\{x\}$ and edge set $E_h \sqcup E_v$.
\item To build $X$, exactly $mn$ squares are attached to $X^{(1)}$, such that
the boundary of each square is of the form $aba'b'$,
where $a,a' \in E_h$, $b,b' \in E_v$. In particular,
the four corners of each square are identified with the vertex $x$.
We denote such a set of $mn$ squares by $R_{m \cdot n}$. 
\item The link $\mathrm{Lk}(X,x)$ of the vertex $x$ in $X$ has to be
isomorphic to the complete bipartite graph on $2m + 2n$ vertices,
where the bipartite structure is induced by the decomposition
of the edges into the two classes $E_h \sqcup E_v$.
Informally speaking, this condition means that
for any $a \in E_h$, $b \in E_v$, the complex $X$ must have a 
unique corner in a unique square
with adjoining edges $a$ and $b$.
\end{itemize}
As a consequence, the universal covering space $\tilde{X}$ of $X$ is the product 
of two regular trees $\mathcal{T}_{2m} \times \mathcal{T}_{2n}$,
see \cite[Proposition~1.1]{BMII} or \cite[Theorem~II.1.10]{Wise}.
By construction, $\Gamma < \mathrm{Aut}(\mathcal{T}_{2m}) \times \mathrm{Aut}(\mathcal{T}_{2n})$
acts freely and transitively on the vertices of $\tilde{X}$,
and for some purposes it is convenient to see $\Gamma$
as a cocompact lattice in 
$\mathrm{Aut}(\mathcal{T}_{2m}) \times \mathrm{Aut}(\mathcal{T}_{2n})$, 
equipped with its usual topology.
Indeed, the main motivation for Burger, Mozes and Zimmer to define
and study such groups $\Gamma$
were expected (super-)rigidity and arithmeticity phenomena
analogous to the famous results for lattices in higher
rank semisimple Lie groups 
(in particular by Margulis \cite{Margulis}).
We will not treat this aspect, but refer to \cite{BMII} and
\cite{BMZ} for interesting developments in this direction.

In the remaining parts of this section we want to discuss
several group theoretic properties of $(2m,2n)$--groups $\Gamma$ needed
in the subsequent sections.

A finite presentation of $\Gamma$ with $m+n$ generators and $mn$ relations
can be directly read off from~$X$:
\begin{align}
\Gamma &= \langle a_1, \ldots , a_m, b_1, \ldots , b_n \mid aba'b' =
1, \text{ for each attached square } aba'b' \rangle \notag \\
&= \langle a_1, \ldots , a_m, b_1, \ldots , b_n \mid R_{m \cdot n}
\rangle \, .\notag
\end{align}
If the $2$-cells of $X$ are metrized as Euclidean squares,
then $X$ is non-positively curved and
$\Gamma$ is a torsion-free CAT(0)-group by \cite[Theorem~4.13(2)]{BrHa}.

Due to the link condition in $X$, 
every element $\gamma \in \Gamma$ can be brought in a unique normal form,
where ``the $a$'s are followed by the $b$'s''.
The idea is to successively replace length $2$ subwords of $\gamma$ of the form $ba$
by $a'b'$, if $a'b'a^{-1}b^{-1} = 1$ in $\Gamma$, or in other words
if (exactly) one of the four squares $a'b'a^{-1}b^{-1}$, $a b'^{-1} a'^{-1} b$,
$a'^{-1} b a b'^{-1}$, $a^{-1} b^{-1} a' b'$ is in $R_{m \cdot n}$.
Analogously, there is a unique normal form, where ``the $b$'s are followed by the $a$'s''.
Here is the precise statement of Bridson-Wise:

\begin{Proposition} \label{normalform}
(Bridson-Wise \cite[Normal Form Lemma~4.3]{BridsonWise})
Let $\gamma$ be any element in a $(2m,2n)$--group
$\Gamma = \langle a_1, \ldots , a_m, b_1, \ldots , b_n \mid R_{m \cdot n} \rangle$.
Then $\gamma$ can be written as 
\[
\gamma = \sigma_a \sigma_b = \sigma_b ' \sigma_a '
\]
where $\sigma_a, \sigma_a '$ are freely reduced words in the subgroup $\langle a_1, \ldots ,a_m \rangle_{\Gamma}$
and $\sigma_b, \sigma_b '$ are freely reduced words in $\langle b_1, \ldots ,b_n \rangle_{\Gamma}$.
The words $\sigma_a, \sigma_a ', \sigma_b, \sigma_b '$ are uniquely determined by $\gamma$.
Moreover, $|\sigma_a| = |\sigma_a '|$ and $|\sigma_b| = |\sigma_b '|$,
where $| \cdot |$ is the word length with respect to the 
symmetric set of standard generators
$\{a_1, \ldots , a_m, b_1, \ldots , b_n\}^{\pm 1}$ of $\Gamma$.
\end{Proposition}
If $\gamma \in \Gamma$ has the form $\sigma_a \sigma_b$ as in 
Proposition~\ref{normalform}, then we say that 
$\gamma$ is in \emph{$ab$-normal form}.
Proposition~\ref{normalform} has some immediate consequences
on the structure of $\Gamma$.

\begin{Corollary} \label{Cor1} 
Let $\Gamma = \langle a_1, \ldots , a_m, b_1, \ldots , b_n \mid R_{m \cdot n} \rangle$ be a $(2m,2n)$--group. Then
\begin{itemize}
\item[$(1)$] The two groups $\langle a_1, \ldots , a_m \rangle_{\Gamma}$ and $\langle  b_1, \ldots , b_n \rangle_{\Gamma}$ 
are free subgroups of $\Gamma$ of rank $m$ and $n$, respectively.
\item[$(2)$] The center $Z \Gamma$ of $\Gamma$ is trivial if $m,n \geq 2$.
\end{itemize}
\end{Corollary}

\begin{proof}
\begin{itemize}
\item[$(1)$] This follows directly from the uniqueness of the normal forms described in Proposition~\ref{normalform}.
\item[$(2)$]
Assume that there is an element $\gamma \in Z \Gamma \setminus \{ 1 \}$ and let 
\[
\gamma = a^{(1)} \ldots a^{(k)} b^{(1)} \ldots b^{(l)}
\]
be its $ab$-normal form,
$a^{(1)}, \ldots, a^{(k)} \in E_h$, $b^{(1)}, \ldots, b^{(l)} \in E_v$, 
where we can assume without loss of generality that  $k \geq 1$ and $l \geq 0$.
Take any element
\[
a \in E_h \setminus \{ a^{(1)}, {a^{(1)}}^{-1} \} \ne \emptyset \, .
\]
(Here, we use $m \geq 2$. Under the assumption $k \geq 0$, $l \geq 1$, we would have used $n \geq 2$.)
Then, we have $a \gamma = \gamma a$, i.e.\
\[
a a^{(1)} \ldots a^{(k)} b^{(1)} \ldots b^{(l)} = a^{(1)} \ldots a^{(k)} b^{(1)} \ldots b^{(l)} a \, .
\]
The left hand side of this equation is already in $ab$-normal form, 
since $a \ne {a^{(1)}}^{-1}$. By uniqueness of the $ab$-normal form, 
we can conclude from the right hand side that
$a = a^{(1)}$, but this is a contradiction to the choice of $a$, and it follows $Z \Gamma = 1$.
\end{itemize}
\end{proof}

For a $(2m,2n)$--group $\Gamma$ we define the homomorphism
$\rho_v : \langle b_1, \ldots, b_n \rangle_{\Gamma} \to
\mathrm{Sym}(E_h)$
as follows. Let $b \in E_v$ and $a \in E_h$, then
$\rho_v (b)(a) := a'$ is the uniquely determined element in
$E_h$ such that $a^{-1} b a' = \tilde{b}$
for some $\tilde{b} \in E_v$.
For a geometric interpretation of $\rho_v$, just draw the
square $a \tilde{b} a'^{-1} b^{-1}$.

Another application of Proposition~\ref{normalform} is the following 
sufficient criterion to show that the centralizer
$Z_{\Gamma}(b) = \{ \gamma \in \Gamma : \gamma b = b \gamma\}$
of $b \in E_v$ is as small as possible.
This will be useful in some results of Section~\ref{SectionCT} and \ref{SectionGamma}.
\begin{Lemma} \label{rhocrit}
Let $\Gamma = \langle a_1, \ldots , a_m, b_1, \ldots , b_n \mid R_{m \cdot n} \rangle$ be a $(2m,2n)$--group. 
Assume that there is an element $b \in E_v$ such that 
$\rho_v (b)(a) \ne a$ for all $a \in E_h$.
Then $Z_{\Gamma}(b) = \langle b \rangle_{\Gamma} \cong \mathbb{Z}$.
\end{Lemma}
\begin{proof}
Obviously, $\langle b \rangle_{\Gamma} < Z_{\Gamma}(b)$. 
We therefore have to show $Z_{\Gamma}(b) < \langle b \rangle_{\Gamma}$. 
Let 
\[
\gamma = a^{(1)} \ldots a^{(k)} b^{(1)} \ldots b^{(l)} \in Z_{\Gamma}(b)
\]
be in $ab$-normal form, 
$a^{(1)}, \ldots, a^{(k)} \in E_h$, $b^{(1)}, \ldots, b^{(l)} \in E_v$, 
$k,l \geq 0$.
Then 
\[
a^{(1)} \ldots a^{(k)} b^{(1)} \ldots b^{(l)} b = b a^{(1)} \ldots a^{(k)} b^{(1)} \ldots b^{(l)} \, .
\]
First assume that $k \geq 1$. The $ab$-normal form of $\gamma b$ starts with $a^{(1)} \ldots a^{(k)}$.
Bringing also $b a^{(1)} \ldots a^{(k)} b^{(1)} \ldots b^{(l)}$ to its
$ab$-normal form, we must have in a first step
$b a^{(1)} = a^{(1)} \tilde{b}$ for some $\tilde{b} \in E_v$, 
i.e.\ ${a^{(1)}}^{-1} b a^{(1)} = \tilde{b} \in E_v$ and therefore
$\rho_v(b)(a^{(1)}) = a^{(1)}$,
which is impossible by assumption, hence $k=0$.
This means $\gamma = b^{(1)} \ldots b^{(l)}$ and
\[
b^{(1)} \ldots b^{(l)} b = b b^{(1)} \ldots b^{(l)} \, .
\]
By uniqueness of the $ab$-normal form of 
\[
b = {b^{(l)}}^{-1} \ldots {b^{(1)}}^{-1}  b b^{(1)} \ldots b^{(l)}
\]
we either have $l = 0$, or $b^{(1)}, \ldots , b^{(l)} \in \{ b, b^{-1} \}$
and hence $\gamma = b^{(1)} \ldots b^{(l)} \in \langle b \rangle_{\Gamma}$.
\end{proof}
Observe that it is very easy to verify for a given set $R_{m \cdot n}$
and $b \in E_v$ if the condition $\rho_v (b)(a) \ne a$ of
Lemma~\ref{rhocrit} holds or not.

We recall the definition of an anti-torus in the context we will use it.
\begin{Definition}
Let $\Gamma = \langle a_1, \ldots , a_m, b_1, \ldots , b_n \mid R_{m \cdot n} \rangle$
be a $(2m,2n)$--group, and $a \in  \langle a_1, \ldots , a_m \rangle_{\Gamma}$,
$b \in  \langle b_1, \ldots , b_n \rangle_{\Gamma}$ two elements.
The subgroup $\langle a, b \rangle_{\Gamma}$ is called an \emph{anti-torus}
in $\Gamma$ if $a,b$ have no commuting non-trivial powers, i.e.\ if 
$a^r b^s \ne b^s a^r$ for all $r, s \in \mathbb{Z} \setminus \{ 0 \}$.
\end{Definition}

\section{Anti-tori in commutative transitive $(2m,2n)$--groups}
\label{SectionCT}
A group $G$ is called \emph{commutative transitive}, if the relation
of commutativity is transitive on the set $G \setminus \{1\}$
(i.e.\ $g_1 g_2 = g_2 g_1$, $g_2 g_3 = g_3 g_2$ always implies
$g_1 g_3 = g_3 g_1$, if $g_1, g_2, g_3 \ne 1$).
Restricting to commutative transitive $(2m,2n)$--groups allows us to give
a very easy criterion to construct anti-tori.
The results stated in this section will be applied to
an interesting subclass of commutative transitive $(2m,2n)$--groups in Section~\ref{SectionGamma}.
\begin{Proposition} \label{Propabantitorus}
Let  $\Gamma = \langle a_1, \ldots, a_m, b_1, \ldots, b_n \mid R_{m \cdot n} \rangle$
be a commutative transitive $(2m,2n)$--group and let
$a \in \langle a_1, \ldots, a_m \rangle_{\Gamma}$,  
$b \in \langle b_1, \ldots, b_n \rangle_{\Gamma}$
be two elements.
Then $\langle a,b \rangle_{\Gamma}$ is an anti-torus in $\Gamma$ if and only if $a$ and $b$ do not commute in $\Gamma$.
\end{Proposition}
\begin{proof}
Assume first that $\langle a,b \rangle_{\Gamma}$ is no anti-torus in $\Gamma$, 
i.e.\ $a^r b^s = b^s a^r$ for some $r,s \in \mathbb{Z} \setminus \{0\}$. 
Obviously, $a$ commutes with $a^r$, and $b$ commutes with $b^s$.
Using the assumption that~$\Gamma$ is commutative transitive, we conclude that $a$ and $b$ commute in $\Gamma$.
The other direction follows immediately from the definition of an anti-torus.
\end{proof}
This gives a dichotomy  for subgroups $\langle a,b \rangle_{\Gamma}$, where $a,b \ne 1$. 
\begin{Corollary} \label{CorZ2orantitorus}
Let $\Gamma = \langle a_1, \ldots, a_m, b_1, \ldots, b_n \mid R_{m \cdot n} \rangle$
be a commutative transitive $(2m,2n)$--group and let
$a \in \langle a_1, \ldots, a_m \rangle_{\Gamma}$,  
$b \in \langle b_1, \ldots, b_n \rangle_{\Gamma}$
be two non-trivial elements.
Then either $\langle a,b \rangle_{\Gamma} \cong \mathbb{Z} \times \mathbb{Z}$ 
or $\langle a,b \rangle_{\Gamma}$ is an anti-torus in $\Gamma$.
\end{Corollary}
\begin{proof}
If $a$ and $b$ do not commute, then $\langle a,b \rangle_{\Gamma}$ is an anti-torus in $\Gamma$ 
by Proposition~\ref{Propabantitorus}.
If $a$ and $b$ commute, then $\langle a,b
\rangle_{\Gamma}$ 
is a finitely generated abelian torsion-free
quotient of $\mathbb{Z} \times \mathbb{Z}$,
hence either $1$, $\mathbb{Z}$ or $\mathbb{Z} \times \mathbb{Z}$.
The first two cases can be excluded by the assumption
$a,b \ne 1$, and using the uniqueness
of the normal forms of powers of $a$ and $b$.
\end{proof}
\begin{Corollary} \label{CorexistAT}
Let $\Gamma = \langle a_1, \ldots, a_m, b_1, \ldots, b_n \mid R_{m \cdot n} \rangle$
be a commutative transitive $(2m,2n)$--group.
Then $\Gamma$ has an anti-torus if and only if $(m,n) \ne (1,1)$.
\end{Corollary}
\begin{proof}
Up to isomorphism, there are only two $(2,2)$--groups: the abelian 
group $\mathbb{Z} \times \mathbb{Z}$, and the (non-commutative transitive) group 
$\langle a_1, b_1 \mid a_1 b_1 a_1 = b_1 \rangle$, where $a_1$
commutes with $b_1^2$. Both groups obviously have no anti-torus.

For the other direction, assume that $(m,n) \ne (1,1)$.
Then there are elements $a \in E_h$ and
$b \in E_v$ which do not commute;
otherwise the $(2m,2n)$--group $\Gamma$ would be a direct product of free groups
\[
\langle a_1, \ldots , a_m \rangle_{\Gamma} \times \langle  b_1, \ldots , b_n  \rangle_{\Gamma} \cong F_m \times F_n \, ,
\]
which is not commutative transitive if $(m,n) \ne (1,1)$.
By Proposition~\ref{Propabantitorus}, $\langle a,b \rangle_{\Gamma}$ is an anti-torus in $\Gamma$.
\end{proof}

The following corollary gives infinitely many anti-tori in $\Gamma$, provided the
centralizer
of some $b$ is cyclic. By Lemma~\ref{rhocrit}, this is for example satisfied for
elements $b \in E_v$ such that 
$\rho_v(b)(a) \ne a$ for all $a \in E_h$.
\begin{Corollary} \label{Corbtantitorus}
Let $\Gamma = \langle a_1, \ldots , a_m, b_1, \ldots , b_n \mid R_{m \cdot n} \rangle$ 
be a commutative transitive $(2m,2n)$--group and let
$b \in \langle b_1, \ldots, b_n \rangle_{\Gamma}$
be an element such that $Z_{\Gamma}(b) = \langle b \rangle_{\Gamma}$. Then $\langle a, b \rangle_{\Gamma}$
is an anti-torus in $\Gamma$ for each $a \in \langle a_1, \ldots, a_m \rangle_{\Gamma} \setminus \{ 1\}$.
\end{Corollary}
\begin{proof}
The assumption $Z_{\Gamma}(b) = \langle b \rangle_{\Gamma}$ implies that $b \ne 1$ and that
$b$ does not commute with any element $a \in \langle a_1, \ldots, a_m \rangle_{\Gamma} \setminus \{ 1\}$.
Now apply Proposition~\ref{Propabantitorus}.
\end{proof}

Similar as for lattices in higher rank semisimple Lie groups,
there is also the important notion of
``reducibility'' and ``irreducibility'' for lattices
acting on a product of trees,
see \cite[Chapter~1]{BMII}:
A lattice $\Gamma < \mathrm{Aut}(\mathcal{T}_{2m}) \times
\mathrm{Aut}(\mathcal{T}_{2n})$
is \emph{reducible} if it is commensurable to a direct product
$\Gamma_1 \times \Gamma_2$ of lattices 
$\Gamma_1 < \mathrm{Aut}(\mathcal{T}_{2m})$,
$\Gamma_2 < \mathrm{Aut}(\mathcal{T}_{2n})$.
Otherwise, $\Gamma$ is called \emph{irreducible}.
Many $(2m,2n)$--groups with interesting group theoretic
properties, like non-residually finite groups
or virtually simple groups (\cite{BMII}, \cite{Rattaggi}), are irreducible,
since reducible $(2m,2n)$--groups contain a subgroup
of finite index which is a direct product of two free groups of finite rank.
There is no known algorithm in general to decide whether
a given $(2m,2n)$--group is irreducible.
However, $(2m,2n)$--groups having an anti-torus
are always irreducible.
\begin{Proposition} (Wise \cite[Section~II.4]{Wise}) \label{PropWiseIrr}
Let $\Gamma = \langle a_1, \ldots , a_m, b_1, \ldots , b_n \mid R_{m \cdot n} \rangle$ be a $(2m,2n)$--group.
If $\Gamma$ has an anti-torus, then it is irreducible.
\end{Proposition}
\begin{proof}
For  $\Gamma < \mathrm{Aut}(\mathcal{T}_{2m}) \times
\mathrm{Aut}(\mathcal{T}_{2n})$
let $\mathrm{pr}_1 : \Gamma \to \mathrm{Aut}(\mathcal{T}_{2m})$
and $\mathrm{pr}_2 : \Gamma \to \mathrm{Aut}(\mathcal{T}_{2n})$
be the two canonical projections. 
Define $\Lambda_1 = \mathrm{pr}_1(\mathrm{ker}(\mathrm{pr}_2)) <
\mathrm{Aut}(\mathcal{T}_{2m})$
and $\Lambda_2 = \mathrm{pr}_2(\mathrm{ker}(\mathrm{pr}_1)) <
\mathrm{Aut}(\mathcal{T}_{2n})$.
Let $\langle a,b \rangle_{\Gamma}$ be an anti-torus in $\Gamma$, 
where $a \in \langle a_1, \ldots , a_m \rangle_{\Gamma}$,
$b \in \langle b_1, \ldots , b_n \rangle_{\Gamma}$,
and suppose that $\Gamma$ is reducible. Then by
\cite[Proposition~1.2]{BMII}, 
the group
$\Lambda_1 \times \Lambda_2$ is a subgroup 
of finite index in $\Gamma$, in particular
the indices
$[\langle a_1, \ldots , a_m \rangle_{\Gamma} : \Lambda_1]$ and
$[\langle b_1, \ldots , b_n \rangle_{\Gamma} : \Lambda_2]$ are finite.
It follows that $a^r \in \Lambda_1$, $b^s \in \Lambda_2$ for some $r,s \in \mathbb{N}$.
But then $a^r b^s = b^s a^r$, a contradiction.
\end{proof}

\begin{Corollary} \label{CorCTirr}
A commutative transitive $(2m,2n)$--group is irreducible if and only if $(m,n) \ne (1,1)$.
\end{Corollary}
\begin{proof}
Any $(2,2)$--group is reducible. If $(m,n) \ne (1,1)$, then 
we combine Corollary~\ref{CorexistAT} and Proposition~\ref{PropWiseIrr}.
\end{proof}

\section{Illustration for the quaternion groups $\Gamma_{p,l}$}
\label{SectionGamma}
For any pair of distinct odd prime numbers $p,l$, we define in this
section a commutative transitive $(p+1,l+1)$--group $\Gamma_{p,l}$,
and can therefore apply the results of Section~\ref{SectionCT}.
With the restriction $p,l \equiv 1 \pmod{4}$, the groups 
$\Gamma_{p,l}$ were originally used by Mozes \cite{Mozes3,Mozes2,Mozes1} to define 
certain tiling systems, so-called two dimensional subshifts of finite
type, and to study a resulting dynamical system.
Later, Burger-Mozes \cite{BMII} used the residually finite group
$\Gamma_{13,17}$ as a building block in the construction of a 
\emph{non}-residually finite $(196,324)$--group and in a construction
of an infinite family of finitely presented torsion-free
virtually simple groups. Kimberley-Robertson \cite{KR}
made explicit computations for many small values of $p,l$,
for example on the abelianization of $\Gamma_{p,l}$.
The condition $p,l \equiv 1 \pmod{4}$ was dropped in \cite{Rattaggi},
and it was shown in \cite{RaRo} that these generalized groups
$\Gamma_{p,l}$ are CSA (i.e.\ all maximal abelian subgroups are
malnormal), in particular they are commutative transitive.

We need some preparation to define the groups $\Gamma_{p,l}$.
For a commutative ring $R$ with unit, let
\[
\mathbb{H}(R) = \{ x_0 + x_1 i + x_2 j + x_3 k : x_0, x_1, x_2, x_3
\in R \}
\]
be the ring of Hamilton quaternions over $R$, i.e.\ $1,i,j,k$ is a free basis,
and the multiplication is determined by $i^2 = j^2 = k^2 = -1$ and
$ij = -ji =k$.
Let $\overline{x} := x_0 - x_1 i - x_2 j - x_3 k \in \mathbb{H}(R)$ be the \emph{conjugate}
of $x =  x_0 + x_1 i + x_2 j + x_3 k \in \mathbb{H}(R)$, and
$|x|^2 := x \overline{x} = \overline{x} x = x_0^2 + x_1^2 + x_2^2 +
x_3^2 \in R$ 
its \emph{norm}.
We write $\Re(x) := x_0$ for the ``real part'' of $x$.

If $R$ is any ring, we denote by $U(R)$ the group of invertible elements 
(with respect to the multiplication) in $R$.

From now on, let $p,l$ be any pair of distinct odd prime numbers.
Let $\mathbb{Q}_p$, $\mathbb{Q}_l$ be the $p$-adic and $l$-adic
numbers, respectively. If $K$ is a field, let as usual 
$\mathrm{PGL}_2(K) = \mathrm{GL}_2(K)/Z\mathrm{GL}_2(K)$,
and write brackets $[A]$ to denote the image of the matrix 
$A \in \mathrm{GL}_2(K)$
under the quotient homomorphism $\mathrm{GL}_2(K) \to \mathrm{PGL}_2(K)$.
We define the homomorphism of groups
\[
\psi_{p,l} : U(\mathbb{H}(\mathbb{Q})) \to 
\mathrm{PGL}_2(\mathbb{Q}_p) \times \mathrm{PGL}_2(\mathbb{Q}_l)
\]
by
\begin{equation}
\begin{split}
\psi_{p,l}(x_0 + x_1 i + x_2 j + x_3 k)=\bigg(&\left[
\begin{pmatrix}
x_0 + x_1 c_p + x_3 d_p & -x_1 d_p + x_2 + x_3 c_p \\
-x_1 d_p - x_2 + x_3 c_p  & x_0 - x_1 c_p - x_3 d_p \\
\end{pmatrix}\right], \notag \\
&\left[\begin{pmatrix}
x_0 + x_1 c_l + x_3 d_l & -x_1 d_l + x_2 + x_3 c_l \\
-x_1 d_l - x_2 + x_3 c_l  & x_0 - x_1 c_l - x_3 d_l \\
\end{pmatrix}\right]
\bigg), \notag
\end{split}
\end{equation}
where $c_p, d_p \in \mathbb{Q}_p$ and $c_l, d_l \in \mathbb{Q}_l$ are elements such that
$c_p^2 + d_p^2 + 1 = 0 \in \mathbb{Q}_p$ and $c_l^2 + d_l^2 + 1 = 0 \in \mathbb{Q}_l$.
This definition is motivated by the following well-known isomorphism:
\begin{Proposition} \label{LemmaDSV}
(see \cite[Proposition~2.5.2]{DSV})
Let $K$ be a field of characteristic different from $2$, and assume that there exist $c, d \in K$ such that
$c^2 + d^2 + 1 = 0$. Then $\mathbb{H}(K)$ is isomorphic to the algebra $M_2(K)$ of $(2 \times 2)$--matrices
over $K$. An isomorphism of algebras is given by the map
\begin{align}
\mathbb{H}(K) &\to M_2(K) \notag \\
x = x_0 + x_1 i + x_2 j + x_3 k &\mapsto 
\left( 
\begin{array}{rr} x_0 + x_1 c + x_3 d & -x_1 d + x_2 + x_3 c \notag \\
-x_1 d - x_2 + x_3 c & x_0 - x_1 c - x_3 d  
\end{array} 
\right) 
\end{align}
and we have
\[
\mathrm{det} \left( 
\begin{array}{rr} x_0 + x_1 c + x_3 d & -x_1 d + x_2 + x_3 c \notag \\
-x_1 d - x_2 + x_3 c & x_0 - x_1 c - x_3 d  
\end{array} 
\right) = |x|^2 \, .
\] 
\end{Proposition}
 
If $p,l \equiv 1 \pmod{4}$, we can choose $d_p = 0$ and $d_l = 0$
in the definition of $\psi_{p,l}$, as in the original definition of Mozes \cite{Mozes3}.
Note that 
\[
U(\mathbb{H}(\mathbb{Q})) = 
\{ x \in \mathbb{H}(\mathbb{Q}) : |x|^2 \in U(\mathbb{Q}) \} =
\mathbb{H}(\mathbb{Q}) \setminus \{0\} \, .
\] 
The homomorphism $\psi_{p,l}$ is not injective, in fact 
\[
\mathrm{ker}(\psi_{p,l}) = ZU(\mathbb{H}(\mathbb{Q})) =
\{x \in U(\mathbb{H}(\mathbb{Q})) : x = \overline{x}\} 
\cong U(\mathbb{Q}) = \mathbb{Q} \setminus \{0\} \, ,
\]
and $\psi_{p,l}(x) = \psi_{p,l}(y)$ if and only if
$y = \lambda x$ for some $\lambda \in U(\mathbb{Q})$.
Observe that 
\[
\psi_{p,l}(x)^{-1} = \psi_{p,l}(x^{-1}) =
\psi_{p,l}\left(\frac{\overline{x}}{|x|^2}\right) =
\psi_{p,l}(\overline{x}) \, .
\]

For an odd prime number $q$, let $X_q$ be the set
\begin{align*}
X_q := \{x = x_0  + &x_1 i + x_2 j + x_3 k \in \mathbb{H}(\mathbb{Z}) \,; 
\quad |x|^2 = q \, ; \notag \\
&x_0 \text{ odd}, x_1, x_2, x_3 \text{ even}, \text{ if } q \equiv 1
\!\!\!\! \pmod 4\,; \notag \\
&x_1 \text{ even}, x_0, x_2, x_3 \text{ odd}, \text{ if } q \equiv 3
\!\!\!\! \pmod 4
\} \, .
\end{align*}
By Jacobi's Theorem (see for example \cite[Theorem~2.1.8]{Lubotzky}), $X_q$ has $2(q+1)$ elements.
Let $Q_{p,l}$ be the subgroup of $U(\mathbb{H}(\mathbb{Q}))$
generated by $X_p \cup X_l \subset \mathbb{H}(\mathbb{Z})$ and $\Gamma_{p,l}$ be its image
$\psi_{p,l}(Q_{p,l})$. Observe that 
\[
\mathrm{ker}(\psi_{p,l}|_{Q_{p,l}}) = \mathrm{ker}(\psi_{p,l}) \cap
Q_{p,l} \cong \{ \pm p^r l^s : r,s \in \mathbb{Z}\} < U(\mathbb{Q}) \, .
\]
Equivalently, $\Gamma_{p,l}$ can be defined as 
\begin{align*}
\psi_{p,l}(\{x \in \mathbb{H}(\mathbb{Z}) \,; 
\quad& |x|^2=p^rl^s, r,s \geq 0\,;  \notag \\
&x_0 \text{ odd}, x_1, x_2, x_3 \text{ even}, \text{ if } |x|^2 \equiv
1 \!\!\!\! \pmod 4\,; \notag \\
&x_1 \text{ even}, x_0, x_2, x_3 \text{ odd}, \text{ if } |x|^2 \equiv
3 \!\!\!\! \pmod 4
\}) \, .
\end{align*}

Note that the set $\psi_{p,l}(X_p)$ has $p+1$ elements, 
since $|X_p| = 2(p+1)$ and $\psi_{p,l}(x) = \psi_{p,l}(-x)$. 
These elements generate a free subgroup 
$\psi_{p,l}(\langle X_p \rangle_{Q_{p,l}}) = \langle a_1, \ldots, a_{\frac{p+1}{2}}\rangle_{\Gamma_{p,l}}$ 
of $\Gamma_{p,l}$ of rank $(p+1)/2$, 
since $\psi_{p,l}(x)^{-1} = \psi_{p,l}(\overline{x})$.
Similarly, $\psi_{p,l}(X_l)$ generates a free subgroup 
$\psi_{p,l}(\langle X_l \rangle_{Q_{p,l}}) = \langle b_1, \ldots, b_{\frac{p+1}{2}} \rangle_{\Gamma_{p,l}}$
of $\Gamma_{p,l}$ of rank $(l+1)/2$.

We summarize the definitions in the following commutative diagram,
where $\psi_{p,l}|$ denotes the restriction of $\psi_{p,l}$ to the respective domain:
\[
\xymatrix{
                        \{ \pm p^r : r \in \mathbb{Z} \} \,   \ar@{^{(}->}[r] \ar@{^{(}->}[d]
&   \, \{ \pm p^r l^s : r,s \in \mathbb{Z} \} \, \ar@{^{(}->}[d]
& \, \{ \pm l^s : s \in \mathbb{Z} \}  \ar@{_{(}->}[l] \ar@{^{(}->}[d] \\
                        \langle X_p \rangle_{Q_{p,l}} \,  \ar@{^{(}->}[r] \ar@{->>}[d]_-{\psi_{p,l}|}    
&   \,  Q_{p,l} \,  \ar@{->>}[d]^-{\psi_{p,l}|}                             
& \, \langle X_l \rangle_{Q_{p,l}}  \ar@{_{(}->}[l] \ar@{->>}[d]^-{\psi_{p,l}|}   \\
                       \langle a_1, \ldots, a_{\frac{p+1}{2}}\rangle_{\Gamma_{p,l}} \, \ar@{^{(}->}[r]
& \, \Gamma_{p,l} \, 
& \, \langle b_1, \ldots, b_{\frac{p+1}{2}} \rangle_{\Gamma_{p,l}} \ar@{_{(}->}[l]
}
\]

Our basic general philosophy is to transfer properties 
of the quaternions to the group $\Gamma_{p,l}$,
and vice versa. For example, the fact that $U(\mathbb{H}(\mathbb{Q}))$
is commutative transitive on non-central elements (cf.\ Lemma~\ref{LeCoQu})
is transferred by Lemma~\ref{LeCo} to the fact (Proposition~\ref{PropCt}) that the group
$\Gamma_{p,l}$ is commutative transitive. 
To simplify the proofs, we introduce the following notation:
If $x = x_0 + x_1 i + x_2 j + x_3 k \in \mathbb{H}(\mathbb{Q})$, let
$\tau (x) := \mathbb{Q} (x_1,x_2,x_3)^T \in \mathbb{Q}^3$.

\begin{Lemma} \label{LeCoQu}
Two quaternions 
$x = x_0  + x_1 i + x_2 j + x_3 k \in \mathbb{H}(\mathbb{Q})$
and $y = y_0  + y_1 i + y_2 j + y_3 k \in \mathbb{H}(\mathbb{Q})$
commute, if and only if $(x_1,x_2,x_3)^T$ and $(y_1,y_2,y_3)^T$
are linearly dependent over $\mathbb{Q}$.
In particular, if $x, y \in \mathbb{H}(\mathbb{Q})$ are two quaternions such that 
$x \ne \overline{x}$ and $y \ne \overline{y}$, then $xy = yx$ if and only if
$\tau(x) = \tau(y)$.
\end{Lemma}
\begin{proof}
The first part follows from the elementary computation
\[
xy - yx =
2(x_2 y_3 - x_3 y_2)i + 2(x_3 y_1 - x_1 y_3)j + 2(x_1 y_2 - x_2 y_1)k \notag \\
= 2
\left| 
\begin{array}{ccc}
i & j & k \\
x_1 & x_2 & x_3 \\
y_1 & y_2 & y_3 
\end{array} 
\right|.
\]
This implies the second part, observing that the condition $x \ne \overline{x}$
is equivalent to the condition $(x_1,x_2,x_3)^T \ne (0,0,0)^T$.
\end{proof}

\begin{Lemma} \label{LeCo}
Two quaternions $x, y \in Q_{p,l}$ commute if and only if 
their images $\psi_{p,l}(x), \psi_{p,l}(y) \in \Gamma_{p,l}$ commute.
\end{Lemma}
\begin{proof}
Obviously $xy = yx$ implies $\psi_{p,l}(x) \psi_{p,l}(y) = \psi_{p,l}(y) \psi_{p,l}(x)$.

For the converse, write as usual 
$x = x_0  + x_1 i + x_2 j + x_3 k \in Q_{p,l} < U(\mathbb{H}(\mathbb{Q}))$,
$y = y_0  + y_1 i + y_2 j + y_3 k \in Q_{p,l}$, and
assume that $\psi_{p,l}(x) \psi_{p,l}(y) = \psi_{p,l}(y) \psi_{p,l}(x)$.
Then $\psi_{p,l}(xy) = \psi_{p,l}(yx)$, hence $xy = \lambda yx$ for some $\lambda \in U(\mathbb{Q})$.
Taking the norm of $xy = \lambda yx$ and using the rule 
$|xy|^2 = |x|^2 |y|^2$, 
we conclude that $1 = |\lambda|^2 = \lambda^2$
and therefore $\lambda = \pm 1$. If $\lambda = -1$, then $xy = - yx$ and
the two general rules $\Re(xy) = x_0 y_0 - x_1 y_1 - x_2 y_2 - x_3 y_3 = \Re(yx)$
and $\Re(-x) = -x_0 = - \Re(x)$ imply that
$\Re(xy) = \Re(-yx) = - \Re(yx) = - \Re(xy)$, hence $\Re(xy) = 0$.
This is impossible, since $\Re(xy) = x_0 y_0 - x_1 y_1 - x_2 y_2 - x_3 y_3$ 
is always an \emph{odd} integer (divided by $p^r l^s$ for some $r, s \geq 0$),
using the parity conditions in the definition of $X_p$ and $X_l$.
Consequently, $\lambda = 1$ and $xy = yx$.
\end{proof}

\begin{Proposition} \label{PropCt}
The group $\Gamma_{p,l}$ is a commutative transitive $(p+1,l+1)$--group.
\end{Proposition} 
\begin{proof}
Mozes showed in \cite[Section~3]{Mozes3} that $\Gamma_{p,l}$ is a
$(p+1,l+1)$--group in the case $p,l \equiv 1 \pmod{4}$. 
It is not difficult to adapt this proof to the general case
(see \cite[Theorem~3.30(5)]{Rattaggi}).

To show that the group $\Gamma_{p,l}$ is commutative transitive,
let $\psi_{p,l}(x), \psi_{p,l}(y), \psi_{p,l}(z) \in \Gamma_{p,l} \setminus \{ 1 \}$
such that $\psi_{p,l}(x)\psi_{p,l}(y) = \psi_{p,l}(y)\psi_{p,l}(x)$ and 
$\psi_{p,l}(y)\psi_{p,l}(z) = \psi_{p,l}(z)\psi_{p,l}(y)$, where
$x,y,z \in Q_{p,l}$ satisfy $x \ne \overline{x}$, $y \ne \overline{y}$ and
$z \ne \overline{z}$.
It follows by Lemma~\ref{LeCo} that $xy = yx$ and $yz = zy$, hence
$\tau(x) = \tau(y) = \tau(z)$ and $xz = zx$ by Lemma~\ref{LeCoQu}.
Thus, we get $\psi_{p,l}(x)\psi_{p,l}(z) = \psi_{p,l}(z)\psi_{p,l}(x)$.
\end{proof}

We can therefore apply the results of Section~\ref{SectionCT} 
to the groups $\Gamma_{p,l}$.
\begin{Proposition} \label{CorGplZf}
Let $\Gamma = \Gamma_{p,l} = \langle a_1, \ldots, a_{\frac{p+1}{2}}, 
b_1, \ldots, b_{\frac{l+1}{2}} \mid R_{\frac{p+1}{2} \cdot \frac{l+1}{2}} \rangle$
and let
$a \in \langle a_1, \ldots, a_{\frac{p+1}{2}} \rangle_{\Gamma}$, 
$b \in \langle b_1, \ldots, b_{\frac{l+1}{2}} \rangle_{\Gamma}$
be two elements. Then
\begin{itemize}
\item[$(1)$]
$\langle a,b \rangle_{\Gamma}$ is an anti-torus in $\Gamma$ 
if and only if $a$ and $b$ do not commute in $\Gamma$.
\item[$(2)$]
If $a,b \ne 1$, then either $\langle a,b \rangle_{\Gamma} \cong \mathbb{Z} \times \mathbb{Z}$ 
or $\langle a,b \rangle_{\Gamma}$ is an anti-torus in $\Gamma$.
\item[$(3)$]
The group $\Gamma$ has an anti-torus and is irreducible.
\item[$(4)$]
If $Z_{\Gamma}(b) = \langle b \rangle_{\Gamma}$, 
then $\langle a, b \rangle_{\Gamma}$ is an anti-torus in $\Gamma$ for each $a \ne 1$.
\end{itemize}
\end{Proposition}
\begin{proof}
\begin{itemize}
\item[$(1)$] Combine Proposition~\ref{Propabantitorus} and Proposition~\ref{PropCt}.
\item[$(2)$] Use Corollary~\ref{CorZ2orantitorus} and Proposition~\ref{PropCt}.
\item[$(3)$] Combine Corollary~\ref{CorexistAT},
  Proposition~\ref{PropCt} and Proposition~\ref{PropWiseIrr}.
\item[$(4)$] This follows from Corollary~\ref{Corbtantitorus} and Proposition~\ref{PropCt}.
\end{itemize}
\end{proof}

\begin{Remark} \label{RemRaRo}
For the group $\Gamma = \Gamma_{p,l}$ there is the following 
easy sufficient criterion 
\cite[Corollary~3.7]{RaRo}
to generate 
elements $b \in \langle b_1, \ldots, b_{\frac{l+1}{2}} \rangle_{\Gamma}$
such that $Z_{\Gamma}(b) = \langle b \rangle_{\Gamma}$.
Take $b \in \langle b_1, \ldots, b_{\frac{l+1}{2}} \rangle \setminus \{1\}$,
write $b = \psi_{p,l}(x_0 + z_0(c_1 i + c_2 j + c_3 k))$,
such that $c_1, c_2, c_3 \in \mathbb{Z}$ are relatively prime, 
$x_0, z_0 \in \mathbb{Z} \setminus \{ 0 \}$, and define $n(b) = c_1^2 + c_2^2 + c_3^2 \in \mathbb{N}$.
If
\[
-\left(\frac{-n(b)}{p}\right) = 1 = \left(\frac{-n(b)}{l}\right) \, ,
\]
then $Z_{\Gamma}(b) = \langle b \rangle_{\Gamma}$, where the 
\emph{Legendre symbol} is
defined as
\begin{equation*} \left(\frac{n}{p}\right)=
\begin{cases}
0 & \text{if $p \mid n$},\\
1 & \text{if $p \nmid n$ and $n$ is a square mod $p$},\\
-1 & \text{if $p \nmid n$ and $n$ is not a square mod $p$}.  
\end{cases}
\end{equation*}
\end{Remark}

Using Lemma~\ref{LeCo}, we also directly get the following result:
\begin{Proposition} \label{Propqat}
Let $x \in \langle X_p \rangle_{Q_{p,l}}$,
$y \in \langle X_l \rangle_{Q_{p,l}}$
be two non-commuting quaternions. 
Then $\langle \psi_{p,l}(x), \psi_{p,l}(y) \rangle_{\Gamma_{p,l}}$ 
is an anti-torus in $\Gamma_{p,l}$.
\end{Proposition}
\begin{proof}
Combine Lemma~\ref{LeCo} and Proposition~\ref{CorGplZf}(1).
\end{proof}

In the other direction, we can use the structure of $\Gamma_{p,l}$
to get a statement on quaternions.
\begin{Proposition}
Let $\Gamma = \Gamma_{p,l} = \langle a_1, \ldots, a_{\frac{p+1}{2}}, 
b_1, \ldots, b_{\frac{l+1}{2}} \mid R_{\frac{p+1}{2} \cdot \frac{l+1}{2}} \rangle$.
Assume that there is an element $b \in E_v$ such that $\rho_v(b)(a) \ne a$ for all $a \in E_h$.
Let $y \in X_l$ such that $\psi_{p,l}(y) = b$. Then there is no
$x \in \langle X_p \rangle_{Q_{p,l}}$ such that $x \ne \overline{x}$ and $xy=yx$. 
\end{Proposition}
\begin{proof}
By Lemma~\ref{rhocrit}, we have $Z_{\Gamma}(b) = \langle b \rangle_{\Gamma}$.
If $x \in \langle X_p \rangle_{Q_{p,l}}$ such that $xy = yx$,
then $\psi_{p,l}(x) \in \langle a_1, \ldots, a_{\frac{p+1}{2}} \rangle_{\Gamma}$
commutes with $\psi_{p,l}(y) = b$, hence
$\psi_{p,l}(x) = 1$ and $x = \overline{x}$.
\end{proof}

We also give an application to number theory in Corollary~\ref{CorMoz}, using the following result of Mozes:
\begin{Proposition} \label{PropMozes315}
(Mozes \cite[Proposition~3.15]{Mozes1})
Let $p,l \equiv 1 \pmod{4}$ be two distinct prime numbers,
\[
\Gamma_{p,l} = \langle a_1, \ldots, a_{\frac{p+1}{2}}, b_1, \ldots, b_{\frac{l+1}{2}} \mid R_{\frac{p+1}{2} \cdot \frac{l+1}{2}} \rangle
\] 
and let $x \in \langle X_l \rangle_{Q_{p,l}}$ such that $x \ne \overline{x}$.
Take $c_1, c_2, c_3 \in \mathbb{Z}$ relatively prime such that
$c := c_1 i + c_2 j + c_3 k \in \mathbb{H}(\mathbb{Z})$ commutes with~$x$.
Then there exists a non-trivial element $a \in \langle a_1, \ldots, a_{\frac{p+1}{2}} \rangle$
commuting with $\psi_{p,l}(x)$ if and only if there are integers $t,u \in \mathbb{Z}$ such that
\[
\mathrm{gcd}(t,u) = \mathrm{gcd}(t,pl) = \mathrm{gcd}(u,pl) = 1
\]
and $t^2 + 4 |c|^2 u^2 \in \{ p^r l^s : r,s \in \mathbb{N} \}$.
\end{Proposition}
\begin{Corollary} \label{CorMoz}
Let $p,l \equiv 1 \pmod{4}$ be two distinct prime numbers
and $\Gamma = \Gamma_{p,l}$.
Let $b = \psi_{p,l}(x_0 + x_1 i + x_2 j + x_3 k) \in \langle b_1, \ldots, b_{\frac{l+1}{2}} \rangle_{\Gamma} \setminus \{ 1 \}$
such that $x_1, x_2, x_3 \in \mathbb{Z}$ are relatively prime.
Moreover, assume that $Z_{\Gamma}(b) = \langle b \rangle_{\Gamma}$.
Then there are no integers $t,u \in \mathbb{Z}$ such that
\[
\mathrm{gcd}(t,u) = \mathrm{gcd}(t,pl) = \mathrm{gcd}(u,pl) = 1
\]
and $t^2 + 4 (x_1^2+x_2^2+x_3^2) u^2 \in \{ p^r l^s : r,s \in \mathbb{N} \}$.
\end{Corollary}

We illustrate some previous results for the group $\Gamma_{5,17}$:
\begin{Corollary}
Let $\Gamma = \Gamma_{5,17}$, $\psi = \psi_{5,17}$ and $b = \psi(3+2i+2j) \in E_v$. Then
\begin{itemize}
\item[$(1)$] The subgroup $\langle \psi(1+2i), \psi(1+4k) \rangle_{\Gamma}$ is an anti-torus in $\Gamma$.
\item[$(2)$] The subgroup $\langle \psi(1+2i), \psi(1+4i) \rangle_{\Gamma} \cong \mathbb{Z} \times \mathbb{Z}$ is no anti-torus in $\Gamma$.
\item[$(3)$] $Z_{\Gamma} (b) = \langle b \rangle_{\Gamma}$.
\item[$(4)$] If $a \in \langle a_1, a_2, a_3 \rangle_{\Gamma} \setminus \{ 1\}$, then 
$\langle a,b \rangle_{\Gamma}$ is an anti-torus in $\Gamma$.
\item[$(5)$] There are no integers $t,u \in \mathbb{Z}$ such that
\[
\mathrm{gcd}(t,u) = \mathrm{gcd}(t,85) = \mathrm{gcd}(u,85) = 1
\]
and $t^2 + 8 u^2 \in \{ 5^r 17^s : r,s \in \mathbb{N} \}$.
\end{itemize}
\end{Corollary}
\begin{proof}
\begin{itemize}
\item[$(1)$] We apply Proposition~\ref{Propqat}, using the obvious fact that $1+2i$ and $1+4k$ do not commute.
\item[$(2)$] The two quaternions $1+2i$ and $1+4i$ commute, 
hence $\psi(1+2i) \in \langle a_1, a_2, a_3 \rangle_{\Gamma} \setminus \{ 1 \}$
commutes with $\psi(1+4i) \in \langle b_1, \ldots, b_9 \rangle_{\Gamma} \setminus \{ 1 \}$.
Now apply Proposition~\ref{CorGplZf}(1),(2).
\item[$(3)$] We check that $\rho_v (b)(a) \ne a$ for all $a \in E_h$ and apply Lemma~\ref{rhocrit}.
Alternatively, we can use \cite[Corollary~3.7]{RaRo} (see Remark~\ref{RemRaRo}), since we have $n(b) = 2$ and 
\[
-\left(\frac{-2}{5}\right) = 1 = \left(\frac{-2}{17}\right) \, .
\]
\item[$(4)$] Apply Proposition~\ref{CorGplZf}(4), using part (3) of this corollary.
\item[$(5)$] This follows from Corollary~\ref{CorMoz}, using part (3) of this corollary and
$b = \psi(3+2i+2j) = \psi(\frac{3}{2}+i+j)$.
\end{itemize}
\end{proof}

Let $b \in \langle b_1, \ldots, b_{\frac{l+1}{2}} \rangle_{\Gamma} \setminus \{ 1 \}$
be a fixed element.
It may happen that $\langle a, b \rangle_{\Gamma_{p,l}}$ is an anti-torus for all $a \in E_h$,
but not for all $a \in \langle a_1, \ldots, a_{\frac{p+1}{2}} \rangle_{\Gamma_{p,l}} \setminus \{ 1\}$.
\begin{Corollary}
Let $\Gamma = \Gamma_{5,7}$, $\psi = \psi_{5,7}$, $b = \psi(1 + 2i + j + k)$, $a_1 = \psi(1 + 2i)$, 
$a_2 = \psi(1 + 2j)$ and $a_3 = \psi(1 + 2k)$.
Then $\langle a, b \rangle_{\Gamma}$ is an anti-torus in $\Gamma$ for all $a \in \{ a_1, a_2, a_3 \}^{\pm 1} = E_h$.
However, $\langle a_2 a_3, b \rangle_{\Gamma}$ is no anti-torus in $\Gamma$.
\end{Corollary}
\begin{proof}
This follows, since $a_2 a_3 = \psi(1+4i+2j+2k)$,
and $1+4i+2j+2k$ commutes with $1 + 2i + j + k$.
\end{proof}

\section{Free anti-tori} \label{SectionFree}
An anti-torus $\langle a, b \rangle_{\Gamma}$ isomorphic to the free
group $F_2$ of rank $2$ is called a \emph{free anti-torus} in $\Gamma$.
It is not known whether there are free anti-tori in $(2m,2n)$--groups,
but we will give in Proposition~\ref{Pfat} a sufficient criterion to
construct free anti-tori in $\Gamma_{p,l}$, using certain free
subgroups
in $U(\mathbb{H}(\mathbb{Q}))$.
An existence theorem for free anti-tori in a class of fundamental groups
of non-positively curved $2$-complexes 
not including $(2m,2n)$--groups, appears in
\cite[Proposition~9.2]{BridsonWise},
but no explicit example of a free anti-torus is given there.
To state our criterion for free anti-tori in $\Gamma_{p,l}$,
we need the following general lemma.
\begin{Lemma} \label{LemmaHom}
Let $\phi : G \to H$ be a homomorphism of groups such that
$\mathrm{ker}(\phi) = ZG$ and let $g_1, \ldots, g_t \in G$, $t \geq 2$. Then
$\langle \phi(g_1), \ldots, \phi(g_t) \rangle_H \cong F_t$ if and only if 
$\langle g_1, \ldots, g_t \rangle_G \cong F_t$.
\end{Lemma}
\begin{proof}
First suppose that $\langle g_1, \ldots, g_t \rangle_G \cong F_t$.
The restriction 
\[
\phi|_{\langle g_1, \ldots, g_t \rangle_G}: \langle g_1, \ldots, g_t \rangle_G \to \langle \phi(g_1), \ldots, \phi(g_t) \rangle_H
\]
is surjective.
It is also injective, since 
\begin{align}
\mathrm{ker}(\phi|_{\langle g_1, \ldots, g_t \rangle_G}) &= \langle g_1, \ldots, g_t \rangle_G \cap \mathrm{ker}(\phi) \notag \\
&= \langle g_1, \ldots, g_t \rangle_G \cap ZG < Z \langle g_1, \ldots, g_t \rangle_G \cong ZF_t = \{ 1 \}, \notag
\end{align}
hence $\langle \phi(g_1), \ldots, \phi(g_t) \rangle_H \cong \langle g_1, \ldots, g_t \rangle_G \cong F_t$.

The other direction is clear, since any relation in $\langle g_1, \ldots, g_t \rangle_G \cong F_t$
induces a relation in $\langle \phi(g_1), \ldots, \phi(g_t) \rangle_H$.
\end{proof}

\begin{Proposition} \label{Pfat}
Let $x \in \langle X_p \rangle_{Q_{p,l}}$, $y \in \langle X_l
\rangle_{Q_{p,l}}$.
Then $\langle \psi_{p,l}(x), \psi_{p,l}(y) \rangle_{\Gamma_{p,l}}$
is a free anti-torus in $\Gamma_{p,l}$ if and only if 
$\langle x,y \rangle_{Q_{p,l}} \cong F_2$.
\end{Proposition}
\begin{proof}
If $\langle x,y \rangle_{Q_{p,l}} \cong F_2$, then
$\langle \psi_{p,l}(x), \psi_{p,l}(y) \rangle_{\Gamma_{p,l}}$
is an anti-torus in $\Gamma_{p,l}$ by Proposition~\ref{Propqat}.
The claim follows now from Lemma~\ref{LemmaHom} applied to the
homomorphism
$\psi_{p,l}|_{Q_{p,l}} : Q_{p,l} \twoheadrightarrow \Gamma_{p,l}$,
where 
\[
\mathrm{ker}(\psi_{p,l}|_{Q_{p,l}}) = \mathrm{ker}(\psi_{p,l}) \cap
Q_{p,l} = ZQ_{p,l} \, .
\]
\end{proof}
We do not know how to apply Proposition~\ref{Pfat} to generate explicit free anti-tori.
Therefore, we pose the following problems:
\begin{Problem}
\begin{itemize}
\item[$(1)$]
Construct a pair $x \in \langle X_p \rangle_{Q_{p,l}}$, $y \in \langle X_l
\rangle_{Q_{p,l}}$ such that
$\langle x,y \rangle_{Q_{p,l}} \cong F_2$.
\item[$(2)$]
Construct a pair $x,y \in \mathbb{H}(\mathbb{Z})$
such that $|x|^2= p^r$, $|y|^2= l^s$ for some $r,s \in \mathbb{N}$
and $\langle x,y \rangle_{U(\mathbb{H}(\mathbb{Q}))} \cong F_2$.
\end{itemize}
\end{Problem}
Nevertheless, we can apply Proposition~\ref{Pfat} in the other
direction
to show that certain $2$-generator groups of quaternions are
\emph{not} free. We first give a general lemma:
\begin{Lemma} \label{Lfi}
Let $\Gamma = \langle a_1, \ldots , a_m, b_1, \ldots , b_n \mid R_{m \cdot n} \rangle$ be a $(2m,2n)$--group
and let $a \in \langle a_1, \ldots , a_m \rangle_{\Gamma}$, $b \in \langle b_1, \ldots , b_n \rangle_{\Gamma}$
be two elements.
If the subgroup $\langle a,b \rangle_{\Gamma}$ has finite index in
$\Gamma$ then $\langle a,b \rangle_{\Gamma} \ncong F_2$.
\end{Lemma} 
\begin{proof}
By \cite{Stallings}, finitely generated, torsion-free, virtually free
groups are free, but $\Gamma$ is clearly not free. 
\end{proof}

This gives an explicit application of Proposition~\ref{Pfat}:
\begin{Proposition} \label{nonfree517}
Let $x = 1 + 2i$ and $y = 1 + 4k$.
Then $\langle x,y \rangle_{U(\mathbb{H}(\mathbb{Q}))} \ncong F_2$.
\end{Proposition}
\begin{proof}
Let $\Gamma = \Gamma_{5,17}$, $a = \psi_{5,17}(x)$ and $b =
\psi_{5,17}(y)$.
Using \textsf{GAP} \cite{GAP}, we check that $\langle a,b \rangle_{\Gamma}$
has index $32$ in $\Gamma$.
By Lemma~\ref{Lfi}, we have $\langle a,b \rangle_{\Gamma} \ncong F_2$,
and Proposition~\ref{Pfat} implies $\langle x,y \rangle_{Q_{5,17}}
\ncong F_2$.
In fact, using the \textsf{GAP}-command
\texttt{PresentationSubgroupMtc}, we have found for example the relation
\begin{align}
&x^3  y^2  x  y^{-1}  x^2  y^{-1} x^2  y^{-1} x^{-4}  y^{-2}  x^{-1}  y  x^{-2}  y^{-1}  x^{-8} y^{-1}  x  y^2  \notag \\
&x y^{-1}  x^{-2} y  x^{-1} y^{-2} x^{-2}  y^{-2}  x^3 y  x^{-2}  y^2  x^2  y^2  x  y^{-1}  x^2  y  x^{-1}  y^{-2} \notag \\ 
&x^{-1}  y  x^8  y  x^2  y^{-1}  x  y^2  x^4  y x^{-2}  y  x^{-2}  y  x^{-1}  y^{-2}  x^{-5}  y^{-1}  x = 1  \notag
\end{align}
of length $106$ in $U(\mathbb{H}(\mathbb{Q}))$.
We do not know if there is a shorter relation.
\end{proof}

We give another example:
\begin{Example}
Let $\psi = \psi_{3,5}$ and 
$\Gamma = \Gamma_{3,5} = \langle a_1, a_2, b_1, b_2, b_3 \mid 
a_1 b_1 a_2 b_2$,  $a_1 b_2 a_2 b_1^{-1}$, 
    $a_1 b_3 a_2^{-1} b_1$, $a_1 b_3^{-1} a_1 b_2^{-1}$,  
    $a_1 b_1^{-1} a_2^{-1} b_3$, $a_2 b_3 a_2 b_2^{-1} \rangle$,
where
\begin{align}
a_1 &= \psi( 1 + j + k), & a_1^{-1} &= \psi( 1 - j - k), \notag \\    
a_2 &= \psi( 1 + j - k), & a_2^{-1} &= \psi( 1 - j + k), \notag \\
b_1 &= \psi( 1 + 2i), & b_1^{-1} &= \psi( 1 - 2i), \notag \\    
b_2 &= \psi( 1 + 2j), & b_2^{-1} &= \psi( 1 - 2j), \notag \\  
b_3 &= \psi( 1 + 2k), & b_3^{-1} &= \psi( 1 - 2k).\notag
\end{align}
Then $\langle a_1, b_1 \rangle_{\Gamma}$ has index $4$ in $\Gamma$
and $\langle a_1^2, b_1^2 \rangle_{\Gamma}$ has index $896$ in
$\Gamma$,
in particular $\langle 1 + j + k, 1 + 2i \rangle_{U(\mathbb{H}(\mathbb{Q}))} \ncong F_2$
and 
\[
\langle (1 + j + k)^2, (1 + 2i)^2
\rangle_{U(\mathbb{H}(\mathbb{Q}))}
=  \langle -1 + 2j + 2k, -3 + 4i
\rangle_{U(\mathbb{H}(\mathbb{Q}))} \ncong F_2 \, .
\]
There is for example the relation
$y  x^3  y^2  x  y^{-1}  x^{-3}  y^{-2}  x^{-1} = 1$, 
where $x = 1 + j + k$ and $y = 1 + 2i$.
\end{Example}

We do not know what happens for increasing powers of $a_1$ and $b_1$:
\begin{Question}
Let $\Gamma = \Gamma_{3,5}$, $a_1 = \psi_{3,5}( 1 + j + k)$ and
$b_1 = \psi_{3,5}( 1 + 2i)$.
\begin{itemize}
\item[$(1)$] Is the index of $\langle a_1^3, b_1^3 \rangle_{\Gamma}$
infinite in $\Gamma$ ?
\item[$(2)$] Is $\langle a_1^3, b_1^3 \rangle_{\Gamma}$ a free
  anti-torus in $\Gamma$ ?
Equivalently, is 
\[
\langle (1 + j + k)^3, (1 + 2i)^3
\rangle_{U(\mathbb{H}(\mathbb{Q}))}
=  \langle -5 + j + k, -11 - 2i 
\rangle_{U(\mathbb{H}(\mathbb{Q}))} \cong F_2 \, ?
\]
\end{itemize}
\end{Question}

There is a more general question of Wise:
\begin{Question}
(\cite[Question 2.7]{Bestvina})
Let $G$ act properly discontinuously and cocompactly on
a CAT(0) space (or let $G$ be automatic). Consider two elements $a$, $b$ of $G$.
Does there exist $n > 0$ such that either the subgroup
$\langle a^n, b^n \rangle_G$ is free or $\langle a^n, b^n \rangle_G$ is abelian?
\end{Question}
Observe that
if $\langle a, b \rangle_G$ is an anti-torus,
then $\langle a^n, b^n \rangle_G$ is never abelian,
and therefore Wise's question in this context is whether there
exists a number $n > 0$ such that $\langle a^n, b^n \rangle_G$
is a free anti-torus.

\section{Free subgroups of $\mathrm{SO}_3(\mathbb{Q})$} \label{SectionSO}
The construction of free subgroups of $\mathrm{SO}_3(\mathbb{R})$
has been studied for example in the context of the Banach-Tarski paradox
(see e.g.\ \cite{Wagon}).
We relate free subgroups of $\mathrm{SO}_3(\mathbb{Q})$ (hence of $\mathrm{SO}_3(\mathbb{R})$)
to free subgroups of $\Gamma_{p,l}$ and to certain free subgroups of $U(\mathbb{H}(\mathbb{Q}))$.

Define $\vartheta: U(\mathbb{H}(\mathbb{Q})) \to \mathrm{SO}_3(\mathbb{Q})$
by mapping $x = x_0 + x_1 i + x_2 j + x_3 k \in U(\mathbb{H}(\mathbb{Q}))$ to the $(3 \times 3)$--matrix 
\[
\frac{1}{|x|^2} \left(
\begin{array}{ccc}
x_0^2 + x_1^2 - x_2^2 - x_3^2  & 2(x_1 x_2 - x_0 x_3)           & 2(x_1 x_3 + x_0 x_2) \\
2(x_1 x_2 + x_0 x_3)           & x_0^2 - x_1^2 + x_2^2 - x_3^2  & 2(x_2 x_3 - x_0 x_1) \\
2(x_1 x_3 - x_0 x_2)           & 2(x_2 x_3 + x_0 x_1)           & x_0^2 - x_1^2 - x_2^2 + x_3^2 \notag
\end{array}
\right) \, .
\]
Note that this is the matrix which represents the $\mathbb{Q}$-linear map $\mathbb{Q}^3 \to \mathbb{Q}^3$,
$y \mapsto xyx^{-1}$ with respect to the standard basis of $\mathbb{Q}^3$,
where the vector $y = (y_1, y_2, y_3)^T \in \mathbb{Q}^3$ is identified with the 
``purely imaginary'' quaternion $y_1 i + y_2 j + y_3 k \in \mathbb{H}(\mathbb{Q})$.
It is well-known that~$\vartheta$ is a surjective homomorphism of groups. Even the restricted map
\[
\vartheta|_{\mathbb{H}(\mathbb{Z}) \setminus \{ 0 \}}: \mathbb{H}(\mathbb{Z}) \setminus \{ 0 \} \to \mathrm{SO}_3(\mathbb{Q})
\]
is surjective, since $\vartheta(\lambda x) = \vartheta(x)$, if $\lambda \in U(\mathbb{Q})$ and $x \in U(\mathbb{H}(\mathbb{Q}))$.
See \cite{LiuRobertson} for an elementary proof of the surjectivity of $\vartheta|_{\mathbb{H}(\mathbb{Z}) \setminus \{ 0 \}}$. 
Moreover, it is easy to check by solving a system of equations that
\[
\mathrm{ker}(\vartheta) = \{ x \in U(\mathbb{H}(\mathbb{Q})) : x = \overline{x} \} = ZU(\mathbb{H}(\mathbb{Q})) \, .
\]
Alternatively, seeing $\vartheta(x)$ as $\mathbb{Q}$-linear map $y \mapsto xyx^{-1}$ as described above, 
we can easily determine the kernel of $\vartheta$ as follows:
\begin{align}
\mathrm{ker}(\vartheta) &= \{ x \in U(\mathbb{H}(\mathbb{Q})) : xyx^{-1} = y, \, \forall y \in \mathbb{H}(\mathbb{Q}) 
\text{ such that } \Re(y) = 0 \} \notag \\
&= \{ x \in U(\mathbb{H}(\mathbb{Q})) : xy = yx, \, \forall y \in \mathbb{H}(\mathbb{Q}) 
\text{ such that } \Re(y) = 0 \} \notag \\
&= \{ x \in U(\mathbb{H}(\mathbb{Q})) : x = \overline{x} \} \cong U(\mathbb{Q})\, . \notag
\end{align}
Observe that if $x \in U(\mathbb{H}(\mathbb{Q})) \setminus ZU(\mathbb{H}(\mathbb{Q}))$, 
then the axis of the rotation $\vartheta(x)$ 
is the line $(x_1, x_2, x_3)^T \cdot \mathbb{Q}$, and the rotation angle $\omega$ satisfies
\[
\mathrm{cos} \, \omega = \frac{x_0^2-x_1^2-x_2^2-x_3^2}{|x|^2} \, ,
\]
or equivalently
\[
\mathrm{cos} \, \frac{\omega}{2} = \frac{x_0}{\sqrt{|x|^2}} \, .
\]

Now, we realize $\Gamma_{p,l}$ as a subgroup of $\mathrm{SO}_3(\mathbb{Q})$,
using the homomorphism $\vartheta$:
\begin{Proposition}
If $\gamma \in \Gamma_{p,l}$, let $x \in Q_{p,l}$ be any quaternion such that
$\psi_{p,l}(x) = \gamma$, and define
$\eta_{p,l}(\gamma) := \vartheta(x)$.
Then $\eta_{p,l}: \Gamma_{p,l} \to \mathrm{SO}_3(\mathbb{Q})$
is an injective homomorphism of groups.
\end{Proposition}
\begin{proof}
We first show that $\eta_{p,l}$ is well-defined, i.e.\ it does not depend
on the choice of $x \in Q_{p,l}$.
Let $x,y \in Q_{p,l}$ such that $\psi_{p,l}(x) = \psi_{p,l}(y) = \gamma$.
Then $y = \lambda x$ for some $\lambda \in U(\mathbb{Q})$,
hence $\vartheta(y) = \vartheta(x)$.

Now we prove that $\eta_{p,l}$ is a homomorphism.
Let $\gamma_1, \gamma_2 \in \Gamma_{p,l}$ 
and $x,y \in Q_{p,l}$ such that $\psi_{p,l}(x) = \gamma_1$,
$\psi_{p,l}(y) = \gamma_2$.
Then $\psi_{p,l}(xy) = \psi_{p,l}(x) \psi_{p,l}(y) = \gamma_1 \gamma_2$
and $\eta_{p,l}(\gamma_1 \gamma_2) = \vartheta(xy) = \vartheta(x) \vartheta(y)
= \eta_{p,l} (\gamma_1) \eta_{p,l} (\gamma_2)$.

Finally, we show that $\eta_{p,l}$ is injective. 
Let $\gamma \in \Gamma_{p,l}$ such that 
$\eta_{p,l}(\gamma) = 1_{\mathrm{SO}_3(\mathbb{Q})}$.
Then $\vartheta(x) = 1_{\mathrm{SO}_3(\mathbb{Q})}$, 
where $x \in Q_{p,l}$ such that $\psi_{p,l}(x) = \gamma$.
It follows that $x \in U(\mathbb{Q})$, hence $\gamma = \psi_{p,l}(x) = 1_{\Gamma_{p,l}}$.
\end{proof}

We therefore have a commutative diagram
\[
\xymatrix{ 
U(\mathbb{Q}) \, \ar@{^{(}->}[d]  & \, \{ \pm p^r l^s : r, s \in \mathbb{Z} \} \ar@{^{(}->}[d] \ar@{_{(}->}[l]  &  & \\
U(\mathbb{H}(\mathbb{Q})) \, \ar@{->}[d]_-{\psi_{p,l}} & 
\, Q_{p,l} \, \ar@{_{(}->}[l] \ar@{^{(}->}[r] \ar@{->}[rd]|-{\vartheta|_{Q_{p,l}}} \ar@{->>}[d]_-{\psi_{p,l}|_{Q_{p,l}}}  &  
\, U(\mathbb{H}(\mathbb{Q})) \ar@{->>}[d]^-{\vartheta} & \\
\mathrm{PGL}_2(\mathbb{Q}_p) \times \mathrm{PGL}_2(\mathbb{Q}_l) &   
\, \Gamma_{p,l} \, \ar@{_{(}->}[l] \ar@{^{(}->}[r]_-{\eta_{p,l}} & \, \mathrm{SO}_3(\mathbb{Q}) \,  \ar@{^{(}->}[r] &  \, \mathrm{SO}_3(\mathbb{R}) \\ 
}
\]

Free subgroups of $Q_{p,l}$, $\Gamma_{p,l}$ and $\mathrm{SO}_3(\mathbb{Q})$ are related as follows:
\begin{Proposition} \label{freeSO}
Let $x^{(1)}, \ldots, x^{(t)}$ be $t \geq 2$ quaternions in $Q_{p,l}$.
Then the following three statements are equivalent
\begin{align}
\langle x^{(1)}, \ldots, x^{(t)} \rangle_{Q_{p,l}} &\cong F_t \\
\langle \psi_{p,l}(x^{(1)}), \ldots, \psi_{p,l}(x^{(t)}) \rangle_{\Gamma_{p,l}} &\cong F_t \\
\langle \vartheta(x^{(1)}), \ldots, \vartheta(x^{(t)})
\rangle_{\mathrm{SO}_3(\mathbb{Q})} &\cong F_t 
\end{align}
\end{Proposition}
\begin{proof}
To show that (1) and (2) are equivalent, we apply Lemma~\ref{LemmaHom}
to the homomorphism $\psi_{p,l}|_{Q_{p,l}} : Q_{p,l} \twoheadrightarrow \Gamma_{p,l}$,
where $\mathrm{ker}(\psi_{p,l}|_{Q_{p,l}}) = ZQ_{p,l}$.

The equivalence between (2) and (3) again follows from Lemma~\ref{LemmaHom}, now applied to the homomorphism
$\eta_{p,l}: \Gamma_{p,l} \to \mathrm{SO}_3(\mathbb{Q})$,
using $\eta_{p,l}(\psi_{p,l}(x)) = \vartheta(x)$ and
$\mathrm{ker}(\eta_{p,l}) = \{ 1 \} = Z\Gamma_{p,l}$.
Note that $Z\Gamma_{p,l} = \{ 1 \}$ holds, since
$\Gamma_{p,l}$ is commutative transitive and non-abelian. In fact,
$Z\Gamma = \{ 1 \}$ holds for any $(2m,2n)$--group $\Gamma$ such that $m,n \geq 2$, 
as seen in Corollary~\ref{Cor1}(2).
\end{proof}

We know some free subgroups of $\Gamma_{p,l}$ and can therefore apply Proposition~\ref{freeSO}.
\begin{Corollary}
Let $\Gamma_{p,l} = \langle a_1, \ldots, a_{\frac{p+1}{2}}, 
b_1, \ldots, b_{\frac{l+1}{2}} \mid R_{\frac{p+1}{2} \cdot \frac{l+1}{2}} \rangle$
and $x^{(1)}, \ldots, x^{(\frac{p+1}{2})} \in X_p$ such that
$\psi_{p,l}(x^{(1)}) = a_1, \ldots, \psi_{p,l}(x^{(\frac{p+1}{2})}) =  a_{\frac{p+1}{2}}$.
Then 
\[
\langle x^{(1)}, \ldots, x^{(\frac{p+1}{2})} \rangle_{Q_{p,l}} \cong F_{\frac{p+1}{2}}
\]
and
\[
\langle \vartheta(x^{(1)}), \ldots, \vartheta(x^{(\frac{p+1}{2})}) \rangle_{\mathrm{SO}_3(\mathbb{Q})} \cong F_{\frac{p+1}{2}} \, .
\]
\end{Corollary}
\begin{proof}
This follows from Proposition~\ref{freeSO}, 
using
\[
\langle a_1, \ldots, a_{\frac{p+1}{2}} \rangle_{\Gamma_{p,l}} \cong F_{\frac{p+1}{2}}
\] 
which holds by Corollary~\ref{Cor1}(1).
\end{proof}
This gives many examples of free groups.
\begin{Example}
Taking the group $\Gamma_{3,5}$, Proposition~\ref{freeSO} implies that
\[
F_2 \cong \langle 1 + j + k, 1 + j - k \rangle_{Q_{3,5}} \, ,
\]
\begin{align}
F_2 &\cong
\langle \vartheta_{3,5}(1 + j + k), \vartheta_{3,5}(1 + j - k) \rangle_{\mathrm{SO}_3(\mathbb{Q})} \notag \\
&= \left\langle
\frac{1}{3}
\left( \begin{array}{rrr}
-1  &  -2            &  2 \\
2  &   1 &  2 \\
-2  &  2 &  1 \notag
\end{array}
\right), \;
\frac{1}{3}
\left( \begin{array}{rrr}
-1  &  2            &  2 \\
-2  &  1 &  -2 \\
-2  &  -2 &  1 \notag
\end{array}
\right) \right\rangle_{\mathrm{SO}_3(\mathbb{Q})} \, ,  \notag
\end{align}
and
\[
F_3 \cong \langle 1 + 2i, 1 + 2j, 1 + 2k \rangle_{Q_{3,5}} \, ,
\]
\begin{gather}
F_3 \cong \langle \vartheta_{3,5}(1+2i), \vartheta_{3,5}(1+2j), \vartheta_{3,5}(1+2k) \rangle_{\mathrm{SO}_3(\mathbb{Q})} =  \notag \\
\left\langle 
\left( \begin{array}{ccc}
1  &  0            &  0 \\
0  &  -3/5 &  -4/5 \\
0  &   \phantom{-}4/5 &  -3/5 \notag
\end{array}
\right), \;
\left( \begin{array}{ccc}
-3/5 & 0 & \phantom{-}4/5 \\
0 & 1 & 0 \\
-4/5 & 0 & -3/5 \notag
\end{array}
\right), \;
\left( \begin{array}{ccc}
-3/5 & -4/5 & 0 \\
\phantom{-}4/5 & -3/5 & 0 \\
0  &  0            &  1 \notag 
\end{array}
\right)
\right\rangle_{\mathrm{SO}_3(\mathbb{Q})} \, .\notag
\end{gather}
\end{Example}
On the other hand, we also get examples of non-free groups:
\begin{Example}
Using Proposition~\ref{freeSO} and Proposition~\ref{nonfree517}, we see that
\[
F_2 \ncong \langle 1 + 2i, 1 + 4k \rangle_{Q_{5,17}}
\]
and
\begin{align}
F_2 &\ncong \langle \vartheta_{5,17}(1 + 2i), \vartheta_{5,17}(1 + 4k) \rangle_{\mathrm{SO}_3(\mathbb{Q})} \notag \\
&= \left\langle
\left( \begin{array}{ccc}
1  &  0            &  0 \\
0  &  -3/5 &  -4/5 \\
0  &   \phantom{-}4/5 &  -3/5 \notag
\end{array}
\right), \;
\left( \begin{array}{ccc}
-15/17  & -\phantom{1}8/17  &  0 \\
\phantom{-1}8/17     & -15/17 &  0 \\
0                &  0              &  1 \notag
\end{array}
\right) \right\rangle_{\mathrm{SO}_3(\mathbb{Q})} \, . \notag
\end{align}
In fact, the long relation in $x^{\pm 1}, y^{\pm 1}$
given in the proof of Proposition~\ref{nonfree517} also holds in $\mathrm{SO}_3(\mathbb{Q})$ for 
the matrices $x = \vartheta_{5,17}(1 + 2i)$, $y = \vartheta_{5,17}(1 + 4k)$.
\end{Example}

\end{document}